\documentclass{amsart}
\usepackage{amscd,amsthm}
\usepackage{latexsym}
\usepackage{capt-of}
\usepackage{graphicx}

\DeclareFontFamily{U}{wncy}{}
\DeclareFontShape{U}{wncy}{m}{n}{<->wncyr10}{}
\DeclareSymbolFont{mcy}{U}{wncy}{m}{n}
\DeclareMathSymbol{\Sha}{\mathord}{mcy}{"58}

\usepackage{amssymb,amsmath}
\usepackage{amsfonts}
\usepackage{comment}
\usepackage{hyperref}
\newcounter{ctfig}

\newcommand{\C}{\mathcal{C}}

\newcommand{\F}{\mathcal{F}}

\newcommand{\JacC}{{\hbox{Jac}_{\lower.5pt\hbox{$_\C$}}}}
\newcommand{\JacF}{{\hbox{Jac}_{\lower.5pt\hbox{$_\F$}}}}

\newcommand{\ord}{\mathop{\rm ord}}

\theoremstyle{plain}

\newtheorem{prop}{Proposition}

\theoremstyle{definition}

\def\F{{\mathbb F}}

\def\C{{\mathbb C}}

\def\e32{{{}_3E_2}}
\def\f32{{{}_3F_2}}
\def\a32{{{}_3A_2}}
\makeatletter
\@namedef{subjclassname@2020}{%
  \textup{2020} Mathematics Subject Classification}
\makeatother

\begin{document}
\bibliographystyle{plain}
\bibstyle{plain}

\title[$Ax^4-By^2=1$]{On an elementary method for solving $Ax^4-By^2=1$.}

\author{P.G. Walsh}
\address{Department of Mathematics\\
University of Ottawa}
\email{gwalsh@uottawa.ca}

\date{\today}

\begin{abstract}
A new method for solving quartic equations due to Luo and Lin is investigated both computationally and theoretically. As a result, a completely straightforward elementary method is given for solving Bumby's equation $3X^4-2Y^2=1$, along with a conjecture, which if resolved, would enable a similar proof for a possibly infinite family of similar equations. 
\end{abstract}

\maketitle

\section{Introduction}
The Diophantine equation $Ax^4-By^2=1$ has a long history, with contributions from many. Most notably, Ljunggren \cite{Lj1} proved many fundamental results on the integer solutions to such equations, while others, such as Cohn \cite{C}, Akhtari \cite{A}, Chen and Voutier \cite{CV}, Bennett, Togb\'e, and the author \cite{BTW} have made noteworthy contributions as well. Much of this topic is covered in the survey \cite{W}.\\

The inspiration for this article is a remarkable result of Richard Bumby \cite{B}.
In 1967, Bumby solved a conjecture of J.H.E. Cohn in \cite{C} by determining that all positive integer solutions to the quartic equation $3x^4-2y^2=1$ are $(x,y)=(1,1)$ and $(x,y)=(3,11)$. The proof given by Bumby uses some clever insights and constructions involving integers in the ring $\mathbb{Z}[\sqrt{-2}]$.\\

We note here, that solving an equation of the form $Ax^4-By^2=1$ can be reduced to the problem of solving one of the form $(t+1)X^4-tY^2=1$, using an argument involving Jacobi symbols given in \cite{Co}. Furthermore, Bumby's result was later generalized in \cite{BTW} using the hypergeometric method to any quartic equation of the form $(t+1)x^4-ty^2=1$, with $t$ of the form $t=u^2+u$, as conjecturally, it is for these specific values of $t$ that such equations have exactly one positive integer solution other than $(x,y)=(1,1)$. In other words, for $t$ not of the form $u^2+u$, it is conjectured that the only positive integer solution is this one obvious solution.\\

Very recently, Lin and Luo \cite{LL} have discovered an elementary method to solve \newline  $(t+1)X^4-tY^2=1$ for $t=1$. Essentially, the method uses a factor base and a modulus to eliminate all but one arithmetic progression determined by the modulus, and then an extraordinary construction to show that the Jacobi symbol of each term in the remaining arithmetic progression evaluated at a certain corresponding integer, always equals $-1$.\\

The purpose of the present paper is to investigate {\em how well} the ideas of \cite{LL} can be used to solve the general equation $(t+1)X^4-tY^2=1$. Though the method of \cite{LL} is exceedingly clever, it appears to have complications and limitations. In particular, for $t=1$, the argument given in \cite{LL} for the evaluation of the relevant Jacobi symbols is somewhat frightful at first glance, with nothing in the argument suggesting that the method has any chance to be generalized. Nevertheless, we will endeavour to develop the ideas. In section 2, we consider the case $t=2$. As we will see, the method works surprisingly well, with the Jacobi symbol argument being far simpler than that in \cite{LL}.

\section{Solving Bumby's equation $3X^4-2Y^2=1$}

\subsection{The Factor Base}
$\; $\\ 

Let $\alpha = \sqrt{3}+\sqrt{2}$, and for $k \ge 1$ define sequences $\{ P_k \}$ and $\{ Q_k \}$ by $$P_k\sqrt{3}+Q_k\sqrt{2} = \alpha^k \; \; \; (k \; {\rm odd})$$
and
$$P_k + Q_k\sqrt{6} = \alpha^k \; \; \; (k \; {\rm even}).$$
\vskip 0.2truecm 
All positive integer solutions to $3x^2-2y^2=1$ in positive integers $(x,y)$ are given by $(x,y)=(P_k,Q_k)$ with $k$ an odd positive integer. The sequence $\{ P_k \}$ can be extended to negative indices by simply taking the corresponding negative power of $\alpha$, and observe that $P_{-k}=P_k$.\\

In order to solve $3x^4-2y^2=1$, we see that this is the same as determining the set of odd positive indices $n$ for which $P_n=x^2$ for some integer $x$. Note that solutions do exist for $n=1$ and $n=3$ (corresponding to $x=1$ and $x=3$ respectively), and hence also for $n=-1$ and $n=-3$. The approach we take is to first exclude a large set of congruence classes for some strategically chosen modulus. These congruence classes are eliminated by primes which serve as witnesses to the non-squareness of $P_n$ with $n$ in a certain congruence class. The modulus is chosen so that the order of the sequence $\{ P_k \}$ modulo each prime in the factor base divides the modulus.\\

Choosing an effective modulus requires some trial and error, which we will refrain from elaborating on, other than to indicate that it must allow for enough primes, yet not grow too large. For the equation under consideration, we have success with the modulus $M=1680$, and the corresponding factor base is given by
$$FB = \{ 11, 13, 29, 41, 43, 59, 71, 83, 89, 97, 109, 179, 211, 241, 337, 419, 587,
673, 881, $$
$\; \; $ $1009, 1901, 3361, 3779, 4549, 5881, 8641, 9601 \}$.\\

A short computation shows that the order of $\{ P_k \}$ modulo each prime in FB divides $1680$, and that, except for the congruence classes 
$$k \in \{ 1,3,837,839,841,843,1677,1679 \}$$ 
modulo $1680$, there is a prime $p \in FB$ for which the Legendre symbol $\large(\frac{P_k}{p}\large) = -1$.\\

We therefore have the following.\\

\newtheorem{Proposition 2.1}{prop}
\begin{prop}
If $P_n=x^2$ for some integer $x$, then $n \equiv \pm 1, \pm 3 \; (\bmod \; 840)$.
\end{prop}

\subsection{The Main Case: $n\equiv 1 \; (\bmod \; 840)$}

We have seen in the previous section that the problem reduces to four congruence classes. Three of these will be handled with relative ease in the final step of the proof. The main thrust of the proof is contained in this section, which consists of the construction of a certain index $b$ related to the index $n$ (in Proposition 1), along with divisibility properties related to the sequences $\{ P_k \}$ and $\{ Q_k \}$, and the evaluation of certain related Jacobi symbols.\\

In order to proceed, we will require a couple of observations regarding the sequences $\{ P_k \}$ and $\{ Q_k \}$. Firstly, if $n$ is odd and $k$ is even, then by comparing the coefficients of $\sqrt{3}$ in
$$\alpha^{n+2k} = P_{n+2k}\sqrt{3}+Q_{n+2k}\sqrt{2}$$ 
with that of
$$\alpha^n \cdot \alpha^{2k} = (P_n\sqrt{3}+Q_n\sqrt{2})\cdot (P_k+Q_k\sqrt{6})^2,$$
we see after expanding the product that 
$$P_{n+2k} \equiv P_n \; (\bmod \; Q_k). \leqno (2.1)$$

The second observation required concerns the quotient $Q_{6k}/Q_{2k}$ with $k$ odd. We will forego the details, and merely state that for $k$ odd
$$Q_{6k}/Q_{2k} = 3(4P_k^2-1)(12P_k^2-1). \leqno (2.2)$$
This is easily seen by raising $\alpha^k$ to the sixth power. The main property to be used in this identity is that for odd $k$, the factor $2P_k+1$ is a divisor of $Q_{6k}$.\\

Assume now that $n \equiv 1 \; (\bmod \; 840)$ ($n \ne 1$), and let $t$ denote the nonzero integer for which $n=1+840t$. This representation will lead to a sequence of four integers $(a,b,c,d)$, depending only on $n$, satisfying very particular properties. It is important to note that $n \in \mathbb{Z}$, and thus may be negative. Define $c=\ord_3(840t)$ and $d=24t/3^c$. Note that the integer $d$ is either $8$ or $16$ modulo $24$. We then define $a,b$ as follows.
$$(a,b) = \begin{cases}
    (35d,3^c) & d \equiv 8 \; (\bmod \; 24) \text{ and $c$ is even}\\
    (7d,5\cdot 3^c) & d \equiv 16 \; (\bmod \; 24) \text{ and $c$ is even}\\
    (5d,7\cdot 3^c) & d \equiv 8 \; (\bmod \; 24) \text{ and $c$ is odd}\\
    (d,35\cdot 3^c) & d \equiv 16 \; (\bmod \; 24) \text{ and $c$ is odd.}
\end{cases}$$
Notice that $a \equiv 16 \equiv -8 \; (\bmod \; 24)$, so we will let $m$ denote an integer for which $a=24m-8$. Notice also that $3|b$, and that $b \equiv 1 \; (\bmod \; 4)$. We will reduce $P_n$ modulo $Q_{6b}$ using (2.1) above as follows.
$$P_n=P_{1+ab}=P_{1+(24m-8)b} \equiv P_{1-8b}=P_{8b-1} \; (\bmod \; Q_{6b}).$$
By (2.2), we have that $2P_b+1$ divides $Q_{6b}$, and so we have shown that
$$P_n \equiv P_{8b-1} \; (\bmod \; (2P_b+1)). \leqno (2.3)$$
The quantity $P_{8b-1}$ can be reduced modulo $2P_b+1$ by first using the identity 
$P_{8b-1}=P_{8b}-2Q_{8b}$, and then expanding $(\alpha^b)^8$ in order to write $P_{8b-1}$ as a polynomial in $P_b$. We forego displaying this calculation, and simply state that after doing so, we obtain the congruence
$$P_n \equiv P_b-2Q_b \; (\bmod \; (2P_b+1)). \leqno (2.4)$$
Since $P_k \equiv 1 \; (\bmod \; 8)$ for all odd $k$, we have that $2P_b+1 \equiv 3 \; (\bmod \; 8)$, giving the following sequence of Jacobi symbol equalities
$$\left( \frac{P_n}{2P_b+1} \right) = \left( \frac{P_b-2Q_b}{2P_b+1} \right) =
-\left( \frac{2P_b-4Q_b}{2P_b+1} \right) = -\left( \frac{-1-4Q_b}{2P_b+1} \right) =
\left( \frac{4Q_b+1}{2P_b+1} \right).$$
Now recall that $b \equiv 1 \; (\bmod \; 4)$, which implies that $b=2r+1$ with $r$ even. In this case, by expanding $\alpha\cdot (\alpha^r)^2$, it is easy to see that $P_{2r+1}=P_{2r}+2Q_{2r}=(2P_r^2-1)+4P_rQ_r$, and that $Q_{2r+1}=P_{2r}+3Q_{2r}=(2P_r^2-1)+6P_rQ_r$.
Therefore,
$$\left(\frac{P_n}{2P_b+1} \right) = \left(\frac{4Q_b+1}{2P_b+1} \right)=\left(\frac{8P_{r}^2+24P_{r}Q_{r}-3}{4P_{r}^2+8P_{r}Q_{r}-1}\right) = \left(\frac{8P_{r}Q_{r}-1}{4P_{r}^2+8P_{r}Q_{r}-1}\right)$$
$$ = -\left(\frac{4P_{r}^2+8P_{r}Q_{r}-1}{8P_{r}Q_{r}-1}\right) = -\left(\frac{4P_{r}^2}{8P_{r}Q_{r}-1}\right) = -1.$$

\subsection{Completion of the Proof}

It remains to deal with the congruence classes $-3,-1$ and $3$ modulo $840$. We first consider the class corresponding to $3$, and then deal with the negative cases to complete the proof.\\

Assume that there is an integer $n\equiv 3 \; (\bmod \; 840)$ ($n \ne 3$) for which $P_n=x^2$. Let $n=3k$. As $P_{3k}=P_k(3P_k^2+6Q_k^2)$, and the only prime that can divide both factors is $3$, there must exist integers $u,v$ for which either $P_k=u^2, 3P_k^2+6Q_k^2=v^2$ or $P_k=3u^2, 3P_k^2+6Q_k^2=3v^2$. The latter case is evidently not possible modulo $8$ after removing the factor of $3$, and so we are left only with the former case. The congruence condition on $n$ forces $k$ to satisfy $k\equiv 1 \; (\bmod \; 280)$. However, by what we have proved so far, the equation $P_k=u^2$ implies that $k \; (\bmod \; 840) \in \{-3,-1,3\}$, which is not possible.\\

If there is an integer $n\equiv -1 \; (\bmod \; 840)$ for which $P_n=x^2$, then because $P_n=P_{-n}$, we can see that $P_n$ is not a square by applying the construction above to $-n$ to produce an integer $b$ with $\left(\frac{P_{-n}}{2P_b+1} \right)=-1$. Finally, the case $n \equiv -3 \; (\bmod \; 840)$ can now be dealt with exactly as we did for the case $n \equiv 3 \; (\bmod \; 840)$.\\

\section{The Extent of the Method}

\subsection{Computational Results Pertaining to the Method}
The result of Akhtari in \cite{A} gives a bound of two positive integer solutions $(x,y)$ to any equation of the form $Ax^4-By^2=1$, which in general is best possible. Despite the existence of the subfamily of such equations with two solutions in \cite{BTW}, most equations conjecturally have only one solution. Therefore, most such quartic equations have not actually been completely solved. In theory then, the method of this paper can be added to the toolkit of ways to solve such equations, over and above all of the existing methods, most of which are implemented in programs such as magma and pari. Moreover, unlike Akhtari's result, which produces an upper bound for all equations, the method here, as we understand it, can only be applied to one equation at a time. In this section, we investigate the degree to which the method is applicable.\\ 

So far, the method has successfully been used to solve $2x^4-y^2=1$ by Lin and Luo \cite{LL}, and in Section 2 above we have shown how the method works well to solve $3x^4-2y^2=1$. Before proceeding to investigate other such equations, we make note of the fact, first proved by Cohn in \cite{Co}, that equations of the form $Ax^4-By^2=1$ with positive integer solutions are in one to one correspondence with equations of the form $(t+1)X^4-tY^2=1$, thus allowing us to focus entirely on this latter family. Computationally, the reduction step involves the computation of the fundamental solution to the pellian-like equation $Ax^2-By^2=1$.\\ 

As described in Section 2, Lin and Luo's method consists of two parts; eliminating almost all arithmetic progressions with witnesses derived from the use of a factor base, and then eliminating the remaining arithmetic progressions with a Jacobi symbol reduction along with subsequent arguments. We begin by describing our experience in applying the first part of the method for equations $(t+1)X^4-tY^2=1$ for $t$ ranging up to some reasonably large bound.\\ 

We remind the reader here that with $\alpha = \sqrt{t+1}+\sqrt{t}$, and for $k \ge 1$, sequences $\{ P_k \}$ and $\{ Q_k \}$ can be defined by 
$$P_k\sqrt{t+1}+Q_k\sqrt{t} = \alpha^k \; \; \; (k \; {\rm odd})$$
and
$$P_k + Q_k\sqrt{t(t+1)} = \alpha^k \; \; \; (k \; {\rm even}).$$

As described in Section 2.1, our approach is to use both a {\em primary} modulus $m$ along with a {\em working} modulus $M$. We typically use $m=840$, and $M=2^r3^sm$ for some non-negative integers $r$ and $s$. For fixed $r$ and $s$, the computation first produces as many primes $p$ as possible with the property that the reduced sequence $P_n \; (\bmod \; p)$ has order which divides $M$. The factor base consists of this set of primes. The second part of the computation computes all indices $j$, with $1 < j < M-1$, for which there is no prime $p$ in the factor base satisfying $\left( \frac{P_j}{p} \right) = -1$. The last part of the computation reduces the indices found in this last part modulo the primary modulus $m$. The expectation of performing this is to have computed the set $\{ 1,839 \}$ when $t$ is not of the form $u^2+u$, and $\{ 1,3,837,839 \}$ otherwise.\\

The result of this computation has been entirely successful in giving the expected result, and almost always with $r$ and $s$ very small. Curiously, it is only for values of $t$ of the form $t=2u^2$ that the computation requires $r$ and $s$ to become somewhat larger, and especially for $t$ being an odd power of $2$.\\

The latter part of the proof deals with the remaining $2$ or $4$ arithmetic progressions, depending on whether $t$ is of the form $t=u^2+u$ or not, and just like in Section 2, the majority of the effort requires dealing with the case $1+840i$ ($i \in \mathbb{Z})$. In what follows, our focus will remain entirely on this case.\\

The method now attempts to show that $P_n=x^2$ is not possible for $n \equiv 1 \; (\bmod \; 840)$ ($n \ne 1$) by using a polynomial $p(x)$, dependent only on $t$, and an input $b$ which depends on $n$, for which the Jacobi symbol $\left( \frac{P_n}{p(b)} \right) = -1$. The choice of $p(x)$ that worked in the aforementioned cases is a factor one of 
$P_{6k}/P_{2k}$ or $Q_{6k}/Q_{2k}$ ($k$ odd) written as polynomials in $P_k$ or $Q_k$ respectively. More precisely, upon computing these polynomial expressions, we determine that the choice of $p(x)$ is a factor of 
$$P(x)=(16(t+1)^2x^4-16(t+1)x^2+1)\cdot (16(t+1)^2x^4-16(t+1)x^2+3).$$
The choice of $b$ described in Section 2 is exactly that used in the work of Lin and Luo \cite{LL}.\\

We performed an extensive computation for $t \le 1000$, using not only the factors of $P(x)$, but many quadratic polynomials, and variations of the described input. Consistently and reliably producing values of the Jacobi symbol all equaling $-1$ only seems to occur for very precise selection of not only $p(x)$ and $b$, but also for $t$. We will elucidate our findings below by way of a conjecture.\\

In short, our computations indicate strongly that the method works only for a very thin set: 
$$t = du^2-1 \; \; {\rm with} \; \;  d \in \{ 2,3,4,6 \}.$$
In particular, with input $b$ as given in Section 2, the following polynomials give the required values of the Jacobi symbol without fail.

\begin{center}
\begin{tabular}{ |c|c|c| } 
 \hline
 $d$ & $t$ & $p(x)$ \\ 
 \hline 
 $2$ & $2i^2-1$ & $8i^2x^2 \pm 4ix -1$ \\ 
 $3$ & $3i^2-1$ & $2ix \pm 1$ \\ 
 $4$ & $4i^2-1$ & $4ix \pm 1$ \\
 $6$ & $6i^2-1$ & $24i^2x^2 \pm 12ix + 1$ \\
 \hline
\end{tabular}
\end{center}

\vskip 0.3truecm 

Proving that the values $p(b)$ are as asserted seems to be very difficult. We are unable so far to prove anything for $t > 2$. We therefore state the following for those wishing to follow up on this endeavour.\\

{\bf Conjecture 3.1}
\vskip 0.1truecm 
For $d=3$, $t=3i^2-1$, $p(x)=2ix+1$ and $i \ge 1$ odd 

$$\left( \frac{P_n}{p(P_b)}\right) = (-1)^{(i-1)/2}\left( \frac{2tQ_b+1}{2iP_b+1} \right)=-1.$$

Note that apart from the case $i=1$, which was proved in Section 2, both equalities given in the conjecture are unproven for all $i \ge 2$. The bottom line is that, should this conjecture be proved, it would enable this elementary proof for an infinite number of equations.

\end{document}